\documentclass[12pt,amsthm]{amsart}

\usepackage{natbib}
\usepackage{lipsum}
\newcommand{\authorfootnotes}{\renewcommand\thefootnote{\@fnsymbol\c@footnote}}%
\makeatother
\usepackage{epsf,latexsym,graphicx}
\usepackage[matrix,arrow,tips,curve]{xy}
\usepackage[leqno]{amsmath}
\usepackage{xy}
\usepackage{amssymb,amsthm}
\usepackage{enumerate}
\usepackage[usenames,dvipsnames]{pstricks}
\usepackage{verbatim}
\usepackage{epsfig}
\usepackage{pst-grad}
\usepackage{pst-plot}
\usepackage{cases}
\newtheorem{definition}{Definition}[section]
\newtheorem{thm}[definition]{Theorem}
\newtheorem{cor}[definition]{Corollary} 
\newtheorem{lem}[definition]{Lemma} 
\newtheorem{prop}[definition]{Proposition} 
 
\newtheorem{example}[definition]{Example}

\newenvironment{pf}{\noindent{}\textbf{Proof.}}{}

\newcommand{\R}{{\ensuremath{\mathbb{R}}}}

\newcommand{\J}{{\ensuremath{\mathcal{J}}}}

\renewcommand{\r}{{\ensuremath{\mathcal{R}}}}

\newcommand\ST{\rule[-1em]{0pt}{2.2em}}

\renewcommand{\Im}{\mathrm{Im\,}}

\date{10 October 2012}

\title[On the dimension of triangular spine spaces]{Homological techniques for the analysis of the dimension of triangular spline spaces}

\author{Bernard Mourrain}
\address{{Bernard Mourrain}
: GALAAD, INRIA Sophia--Antipolis M\'editerran\'ee, France}
\email{Bernard.Mourrai@inria.fr}

\author{Nelly Villamizar}
\address{{Nelly Villamizar}: Centre of Mathematics for Applications, University of Oslo, Norway}
\email{nellyyv@cma.uio.no}

\begin{document}
\maketitle
\begin{abstract}
The spline space $C_k^r(\Delta)$ attached to a subdivided domain
$\Delta$ of $\R^{d} $ is the vector space of functions of class
$C^{r}$ which are polynomials of degree $\le k$ on each piece of
this subdivision.  Classical splines on planar rectangular grids
play an important role in Computer Aided Geometric Design, and spline
spaces over arbitrary subdivisions of planar domains are now
considered for isogeometric analysis applications.  We address the
problem of determining the dimension of the space of bivariate
splines $C_k^r(\Delta)$ for a triangulated region $\Delta$ in the
plane. Using the homological introduced by \cite{bil}, we
number the vertices and establish a formula for
an upper bound on the dimension. There is no restriction on the
ordering and we obtain more accurate approximations to the
dimension than previous methods and furthermore, in certain cases even an exact value can be found. The construction makes also possible to get a short proof
for the dimension formula when $k\ge 4r+1$, and the same method we
use in this proof yields the dimension straightaway for many
other cases.
\end{abstract}

\section{Introduction}

Let $\Delta$ be a connected, finite two dimensional simplicial
complex, supported on $|\Delta|\subset\R^2$, with $|\Delta|$ homotopy
equivalent to a disk.  We denote by $C_k^r(\Delta)$ the vector space
of all $C^r$ functions on $\Delta$ that, restricted to any simplex in
$\Delta$, are given by polynomials of degree less or equal than
$k$. These functions are called splines, and they have many practical
applications, including the finite element method for solving
differential equations (\cite{strangfix}). Recently they have also been considered for
isogeometric analysis applications (\cite{iga-chb-09}), in
computer-aided design for modeling surfaces of arbitrary topology (\cite{f-cscag-93}). In
these application areas, spline functions of low degree are of
particular interest, and the degree of smoothness attainable is an
important design consideration.

A fundamental problem is to determine the dimension of this vector
space as a function of known information about the subdivision
$\Delta$.  A serious difficulty for solving this problem is that the
dimension of the space $C_k^r(\Delta)$ can depend not only on the
combinatorics of the subdivision, but also on the geometry of
the triangulation i.e., how $|\Delta|$ is embedded in $\R^{2}$.  In
\cite{Sch}, the author presented a lower and an upper bound on the
dimension of spline spaces of arbitrary degree and smoothness for
general triangulations;  using Bernstein-B\'ezier methods, a formula for the dimension for
$k\geq 4r+1$  was obtained by \cite{alfschu}. The result was extended to $k\geq 3r+2$
in \cite{hd}, and in \cite{alfschu3r+1} a dimension formula is proved for almost all triangulations for $k=3r+1$. However, there are no explicit formulas for the
dimension of the spline spaces $C_k^r(\Delta)$ for degree $k<3r+2$ for
general triangulations.  The use of homological algebra in spline
theory was introduced by \cite{bil}. He obtained the dimension
of $C_k^1(\Delta)$ for all $k$ for generic triangulations.  
\cite{Local}, introduced a chain complex different from that used by Billera; this complex was studied by \cite{Local}, \cite{family}, \cite{Sck}, and \cite{gs}.
The lower homology modules of the chain complex in this construction,
differ from the one introduced by  Billera, and they have nicer
properties.  The connection between fat points and the spline space
defined on $\Delta$, allows to give a complete characterization of the
free resolutions of the ideals generated by power of linear forms
appearing in the chain complex, and so to prove the dimension formula
for $C_k^r(\Delta)$ for sufficiently high degree (\cite{family}). 

The main contribution of the paper is a new formula for an upper bound on the
dimension of the spline space. The formula applies to any ordering 
established on the interior vertices of the partition, contrarily to the upper
bound formulas in \cite{Sch}, \cite{LaiSchu}. Having no restriction on
the ordering makes it possible to obtain accurate approximation to
the dimension and even exact value in many cases.  As a consequence,
we give a simple proof for dimension formula 
when $k\geq 4r+1$. 

The paper is structured as follows.
In Section \ref{construction} we recall the construction of this chain complex
and some of the properties of the homology modules. 
We describe the dimension as a sum of a combinatorial part and the
dimension of an homology module.
This latter term which happens to be zero in high degree is always non-negative, so the
formula gives a lower bound for the dimension of the spline space for
any degree $k$.
We present the main result in Section
\ref{sectionUpperBound}, which is an upper bound on the
dimension of the spline space. 
In Section \ref{comp} we
compare the formulas of the lower and the upper bound with those
appearing in \cite{Sch}. 
The result about the exact dimension
for $k\geq 4r+1$ is proved in Section
\ref{exactdimension}. This latter result and some other examples that we present in
Section \ref{examples}, illustrate the interest of the homology
construction for proving exact dimension formulas. 

\section{Construction of the chain complex}\label{construction}

We reproduce some notations and definitions presented in \cite{Local}, restricting them  to the case where $\Delta$ is a planar simplicial complex supported on a disk.

Denote by $\Delta^0$ the set of interior faces of $\Delta$, and by
$\Delta_{i}^0$ ($i=0,1,2$) the set of $i$-dimensional interior faces
whose support is not contained in the boundary of $|\Delta|$.  We
denote by $f_i^0$ the cardinality of these sets, and by
$\partial\Delta$ the complex consisting of all $1$-faces lying on just
one $2$-face as well as all subsets of them.  As it will be convenient
to study the dimension of the vector space $C_k^r(\Delta)$, we embed
$\Delta$ in the plane $\{z=1\}\subseteq\R^3$ and form the cone
$\hat{\Delta}$ over $\Delta$ with vertex at the origin. Denote by
$C_k^r(\hat{\Delta})$ the set of splines on $\hat{\Delta}$ of
$C^r$-smoothness and degree exactly $k$. Then
$C^r(\hat{\Delta}):=\oplus_{\geq 0}C_k^r(\hat{\Delta})$ is a graded
$\R$-algebra and there is an isomorphism of $\R$-vector spaces between
$C_k^r(\Delta)$ and the elements in $C^{r}(\hat{\Delta})$ of degree
exactly $k$ (\cite{BR1}), in particular
\[\dim C_k^r(\Delta)=\dim C_k^r(\hat{\Delta}).\]
Define $R := \R[x, y, z]$. For an edge $\tau\in\Delta_1^0$, let
$\ell_\tau$ denote a non-zero homogeneous linear form vanishing on
$\hat\tau$, and define the ideal $\J(\beta)$ of $R$ for each simplex
$\beta\in\Delta^0$ as follows:
\begin{align*}
\J(\sigma)&=\langle 0\rangle&&\text{for each} \;\sigma\in\Delta_2^0\\
\J(\tau)&=\langle \ell_{\tau}^{r+1}\rangle &&\text{for each} \;\tau\in\Delta_1^0\\
\J(\gamma)&=\langle \ell_{\tau}^{r+1}\rangle_{\tau\ni\gamma}&&\text{for each} \;\gamma\in\Delta_0^0, \,\tau\in\Delta_{1}^0.
\end{align*}
Consider the chain complex $\mathcal{R}$ defined on $\Delta^0$  as $\mathcal{R}_i=R^{f_{i}^{0}}$, with the  usual simplicial boundary maps $\overline{\partial_i}$ used to compute the relative (modulo $\partial\Delta)$ homology with coefficients in $R$ as in \cite{bil}.
Let $\r/\J$ be the chain complex obtained as the quotient of $\r$ by $\J$,
\vspace{0.1cm}
\[0\longrightarrow\bigoplus_{\sigma\in\Delta_2}\mathcal{R}\xrightarrow{\partial_2}\bigoplus_{{\tau\in\Delta}_1^{0}}\mathcal{R/J}(\tau)\xrightarrow{\partial_1}\bigoplus_{{\gamma\in\Delta}_0^{0}}\mathcal{R/J}(\gamma)\longrightarrow 0\]
\vspace{0.1cm}where the maps $\partial_i$ are the induced by the simplicial boundary maps $\overline{\partial_i}$.
The complex $\r/\J$ was introduced by \cite{Local}, and agrees
with the complex studied by Billera except at the vertices.
It was shown in \cite{bil} that $C^r(\hat{\Delta})$ is isomorphic to the top homology module of $\r/\J$
\[H_{2}(\mathcal{R/J}):=\ker(\partial_2).\] Let us consider the short exact sequence of complexes
\[0\longrightarrow\mathcal{J}\longrightarrow\mathcal{R}\longrightarrow\mathcal{R/J}\longrightarrow 0\]
that gives rise to the long exact sequence of homology modules
\begin{align}\label{longsequence}
0\rightarrow H_2(\mathcal{R})\rightarrow &H_2(\mathcal{R}/\J)\rightarrow H_1 (\J)\rightarrow H_1(\mathcal{R})\nonumber\\
&\rightarrow H_1(\mathcal{R}/\J)\rightarrow H_0(\J)\rightarrow H_0(\mathcal{R})\rightarrow H_0 (\mathcal{R}/\J)\rightarrow 0
\end{align}
Since $\Delta$ is supported on a disk, both $H_0(\r)$ and $H_1(\r)$ are zero. Hence, by the long exact sequence, it implies that $H_{0}(\mathcal{R/J})$ is also zero  and $H_{1}(\mathcal{R/J})$ is isomorphic to $H_{0}(\mathcal{J})$.
Applying the Euler characteristic equation to the chain complex $\r/\J$
\[\chi(H(\mathcal{R/J}))=\chi(\mathcal{R/J}),\]
and considering the modules in degree exactly $k$, leads to the formula
\begin{equation}\label{eq}
\dim C^r_k(\Delta)=\sum_{i=0}^2 (-1)^i\hspace{-0.1cm}
\sum_{\beta\in\Delta_{2-i}^0}\dim \mathcal{R}/\mathcal{J}(\beta)_k + \dim H_{0}(\J)_k.
\end{equation}
The aim is to determine the modules in the previous formula as functions of known information about the subdivision $\Delta$.

We know that
\begin{align}
&\sum_{\sigma\in\Delta_{2}^{0}}\dim  \mathcal{R}_k=f_2^0 \,\binom{k+2}{2}\label{triang}\\
&\sum_{\tau\in\Delta_1^{0}}\dim \mathcal{R}/ \mathcal{J}(\tau)_k=f_1^{0}\,\biggl[\binom{k+2}{2}-\binom{k+2-(r+1)}{2}\biggr].\label{edges}
\end{align}
For the computation of  $\dim \mathcal{R}/\mathcal{J}(\gamma)_k$, \cite{family} proposed the following resolution for $\mathcal{R}/\mathcal{J}(\gamma)_k$.  Without loss of generality we translate $\gamma$ to the origin and assume that the linear forms in $\mathcal{J}(\gamma)$ involve only variables $x,y$. Letting $\ell_1^{r+1},\dots ,\ell_{t}^{r+1}$ be a minimal
generating set for $\mathcal{J}(\gamma)$,  then a free resolution for $\mathcal{R}/\mathcal{J}(\gamma)_k$ is given by:
{\small
\begin{equation}\label{resolution}
0\rightarrow \mathcal{R}(-\Omega-1)^{a}\oplus \mathcal{R}(-\Omega)^{b}\rightarrow\oplus_{j=1}^{t}\mathcal{R}(-r-1){\longrightarrow}
\mathcal{R}\rightarrow \mathcal{R}/\mathcal{J}(\gamma)\rightarrow 0
\end{equation}}where $\Omega-1$ is the socle degree of $\mathcal{R}/\mathcal{J}(\gamma)$; $\Omega$ and the multiplicities $a$ and $b$ are given by
\begin{equation}\label{formulas}
\Omega = \left\lfloor \frac{t\,r}{t-1}\right\rfloor + 1,\quad a=t\,(r+1)+(1-t)\,\Omega,\quad b=t-1-a.
\end{equation}
In the case $t=1$, we take $A=b=\Omega=0$ so that $\mathcal{R}(-\Omega-1)^{a}\oplus \mathcal{R}(-\Omega)^{b}=0$.
Applying this to each vertex $\gamma_i\in\Delta_0^0$, and letting $t_i$, $\Omega_i$, $a_i$ and $b_i$ denote the values for $t$, $\Omega$, $a$ and $b$ at $\gamma_i$ respectively,  leads to
\begin{align}\label{vertices}
\dim\bigoplus_{\gamma_i\in\Delta_0^0}\mathcal{R/J}(\gamma_i)_k=\sum_{i=1}^{f_0^0}\biggl[&\binom{k+2}{2}-t_i\binom{k+2-(r+1)}{2}+\nonumber\\
&b_i\,\binom{k+2-\Omega_i}{2}+a_i\,\binom{k+2-(\Omega_i+1)}{2}\biggr].
\end{align}
Let us notice that the previous formula does not change when we take $t_i$ as the number of different slopes of the edges containing the vertex $\gamma_i$. Then
 (\ref{triang}), (\ref{edges}), and (\ref{vertices}), and the fact that $\dim H_1(\r/\J)_k\geq 0$, yields the following theorem:
\begin{thm} \label{lowbHOM}
The dimension of $C^{r}_{k}(\Delta)$ is bounded below by
{\small
\begin{align*}
&\dim C_k^r(\Delta)\geq
\binom{k+2}{2}+f_1^{0}\,\binom{k+2-(r+1)}{2}\\
&-\sum_{i=1}^{f_0^0}\biggl[t_i\binom{k+2-(r+1)}{2}-b_i\,\binom{k+2-\Omega_i}{2}-a_i\,\binom{k+2-(\Omega_i+1)}{2}\biggr],
\end{align*}}where $t_{i}$ is the number of different slopes of the edges containing the vertex $\gamma_{i}$, and
\[\Omega_i=\biggl\lfloor\frac{t_i\,
  r}{t_i-1}\biggr\rfloor+1,\quad a_i=t_i\, (r+1)  +(1-t_i)\,\Omega_i,\quad b_i=t_i-1-a_i.\]
\end{thm}
In \cite{Sck}, it was proved that  the homology module $H_1(\mathcal{R/J})_k$  vanishes for sufficiently high degree, thus the lower bound in the latter theorem is actually the exact dimension formula for $C_k^r(\Delta)$ when $k\gg 0$  \cite[Theorem 4.2]{gs}.

\section{An upper bound on the dimension of $C_k^r(\Delta)$}\label{sectionUpperBound}
Let us fix an ordering
$\gamma_{1}, \ldots, \gamma_{f_{0}^{0}}$
for the vertices in $\Delta^0_0$. For each vertex $\gamma_i$, denote by ${N} (\gamma_i)$ the set edges that contain this vertex, and  define
$\tilde{t}_i$ as the number of different slopes of the edges connecting $\gamma_i$ to one of the first $i-1$ vertices in the list or to a vertex on the boundary.

\begin{thm}\label{upbHOM}
The dimension of $C^{r}_{k}(\Delta)$ is bounded by
{\small
\begin{align*}
\dim \; &C_k^r(\Delta)\leq
\binom{k+2}{2}+f_1^{0}\,\binom{k+2-(r+1)}{2}\\
&\hspace{-0.7cm}-\sum_{i}\biggl[\tilde{t}_i\binom{k+2-(r+1)}{2}-\tilde{b}_i\,\binom{k+2-\tilde{\Omega}_i}{2}-\tilde{a}_i\,\binom{k+2-(\tilde{\Omega}_i+1)}{2}\biggr]
\end{align*}}with $\tilde{t}_i$ as we have defined above and
\[\tilde\Omega_i=\biggl\lfloor\frac{\tilde{t}_i\, r}{\tilde{t}_i-1}\biggr\rfloor+1,\quad
\tilde{a}_i=\tilde{t}_i\, (r+1)
+(1-\tilde{t}_i)\,\tilde\Omega_i,\quad
\tilde{b}_i=\tilde{t}_i-1-\tilde{a}_i.\] if $\tilde{t}_{i}>1$ and
$\tilde{a}_{i}=\tilde{b}_{i}=\tilde{\Omega}_{i}=0$ if $\tilde{t}_{i}=1$ .
\end{thm}

\begin{pf}
By the long exact sequence (\ref{longsequence}), and the fact that $H_1(\r)=0$, we have the short exact sequence
\[0\longrightarrow H_2(\r)\longrightarrow H_2(\r/\J)\longrightarrow H_1(\J)\longrightarrow 0\]
The Euler characteristic equation applied to the complex $\r$ implies that $H_2(\r)_k=R_k$. Hence the isomorphism between $C^r(\hat\Delta)$ and the homology module $H_2(\r/\J)$ implies that
\[\dim C_k^r(\Delta)=\dim R_k+\dim H_1(\J)_k\]
where the complex of ideals $\J$ (as defined above) is given by
\[0\longrightarrow\bigoplus_{\tau\in\Delta_1^0}\J(\tau)\xrightarrow{\partial_1}\bigoplus_{\gamma\in\Delta_0^0}\J(\gamma)\longrightarrow 0\]
and where $H_{1} (\J) = \ker \partial_{1} = K_{1}$.
Define $W_1:=\Im(\partial_1)$. By the exact sequence
\[0\longrightarrow K_1 \longrightarrow \bigoplus_{\tau\in\Delta_1^0}\J(\tau)\longrightarrow W_1\longrightarrow 0\]
we get \[\dim  C_k^r(\Delta)= \dim R_{k} + \sum_{\tau\in\Delta_1^0}\dim \J(\tau)_k-\dim (W_{1})_k\]
Therefore, to find an upper bound on $\dim C_k^r(\Delta)$ it is enough
to find a lower bound on the dimension of $W_{1}= \Im \partial_{1}$ in
degree $k$.
Let us consider the map
\[
\delta_1:\bigoplus_{\tau= (\gamma,\gamma') \in\Delta_1^0} \J(\tau)\, [\tau]
\to
\bigoplus_{\gamma_i\in\Delta_0^0}
\bigoplus_{\tau \in N (\gamma)}
R\, [\tau|\gamma]
 \]
such that $\delta_{1} ([\tau]) = [\tau|\gamma] -[\tau|\gamma']$ for
$\tau= (\gamma,\gamma') \in \Delta_{1}^{0}$,
and the map
\[\varphi_1:
\bigoplus_{\gamma_i\in\Delta_0^0}
\bigoplus_{\tau \in N (\gamma)}
R\, [\tau|\gamma]
\to
\bigoplus_{\gamma\in\Delta_0^0}
R\, [\gamma]  \]
with
\begin{align*}
{\varphi}_{1} ([\tau| \gamma]) &= [\gamma] \ \mathrm{if}\ \gamma\in \Delta^{0}_{0},\\
&= 0 \ \;\; \mathrm{if}\  \gamma\not\in \Delta^{0}_{0}.
\end{align*}
Then, we have $\partial_{1} = \varphi_{1} \circ \delta_{1}$.
We consider now the map
\[\pi_1:
\bigoplus_{\gamma\in\Delta_0^0}
\bigoplus_{\tau \in N (\gamma)}
R\, [\tau|\gamma]
\to
\bigoplus_{\gamma\in\Delta_0^0}
\bigoplus_{\tau \in N (\gamma)}
R\, [\tau|\gamma]
\]
with
${\pi}_{1} ([\tau| \gamma]) = 0$ if $\gamma$ is the end
point of biggest index of $\tau$
and ${\pi}_{1} ([\tau| \gamma]) = [\tau| \gamma]$ otherwise.
We denote by $\tilde{\partial}_{1} = \varphi_{1} \circ \pi_{1} \circ \delta_{1}$.

For $\gamma\in \Delta_{0}^{0}$, we denote by $\tilde{N} (\gamma)$ the set of interior edges $\tau$ connecting $\gamma$
to another vertex which is not of bigger index. Let $\tilde{\J} (\gamma) = \sum_{\tau\in \tilde{N} (\gamma)} R \, \ell_{\tau}^{r+1} \subset \J (\gamma)$.
By construction, we have
$$
\Im \tilde{\partial}_{1} = \bigoplus_{\gamma\in\Delta_0^0} \widetilde{\J}(\gamma)\, [\gamma]
$$
and $\dim (W_{1})_{k} = \dim (\Im \partial_{1})_{k} \ge \dim (\Im \tilde{\partial}_{1})_{k}.$

Using formulas \eqref{triang}, \eqref{edges} and the resolution
\eqref{resolution} applied to $\tilde{\J} (\gamma)$, we obtain the
upper bound
\begin{align*}
&\dim C_k^r(\Delta)\leq
\binom{k+2}{2}+f_1^{0}\,\binom{k+2-(r+1)}{2}\\
&-\sum_{i=1}^{f_0^0}\biggl[\tilde{t}_i\binom{k+2-(r+1)}{2}-\tilde{b}_i\,\binom{k+2-\tilde{\Omega}_i}{2}-\tilde{a}_i\,\binom{k+2-(\tilde{\Omega}_i+1)}{2}\biggr],
\end{align*}
with $\tilde{t}_{i} = |\tilde{N} (\gamma_{i})|$ and $\tilde{\Omega}_{i},
\tilde{a}_{i}, \tilde{b}_{i}$ defined as in \eqref{formulas}.
\end{pf}

As an immediate consequence of this theorem we mention the following  result.
\begin{cor}\label{>r+1}
If for a numbering on the vertices in $\Delta_0^0$, either $t_i=\tilde t_i$ or $\tilde t_i\geq {r+1}$  for every vertex $\gamma_i$, then the upper bound we get, equals the lower bound, obtaining so the exact dimension formula for the spline space.
\end{cor}
\begin{pf}
  We compare the terms corresponding to each interior vertex in both
  formulas. If $t_i=\tilde t_i$ then they are trivially the same.  By
  definition (\ref{formulas}), for $t\geq r+1$, $\Omega$ and $a$ are
  both constant equal to $r+1$. Hence if $r+1\leq\tilde t_i$ $(\leq
  t_i)$, the terms in the binomials are the same when we compare them
  in the formulas in Theorem \ref{lowbHOM} and \ref{upbHOM}. Since the
  respective terms in $t_i$ and $\tilde t_i$ cancel out then we get
  the equality of the bounds.
\end{pf}
\section{The bounds on $\dim C_k^r(\Delta)$ given by \cite{Sch}}\label{comp}

In this section we compare the bounds on $\dim C_k^r(\Delta)$ found by \cite{Sch}, with the lower and upper bounds  given in the previous two sections.
With the  notation as before, the upper bound presented by Schumaker can be stated as follows.

\begin{thm}\cite[Theorem 2.1]{Sch}\label{SchTheorem}
Suppose that the vertices $\gamma_i\in\Delta^0_ 0$ of the partition are numbered in such a way that each pair of consecutive vertices in the list are corners of a common triangle in $\Delta$. For each  $\gamma_i$ define $\tilde{t}_i$ as the number of edges with different slopes joining the vertex $\gamma_i$ to a vertex in the boundary of $\Delta$ or to one of the first $i-1$ vertices. Then
{\small
\begin{align}\label{upperboundSchumaker}
\dim C_k^r(\Delta)&\leq\binom{k+2}{2}+ f_1^0\binom{k-r+1}{2} \nonumber\\
&-f_0^0\biggl[\binom{k+2}{2}-\binom{r+2}{2}\biggr]+\sum_{i=1}^{f_0^0}\hspace{0.1cm}\sum_{j=1}^{k-r}(r+j+1-j\cdot\tilde t_i)_+.
\end{align}}
\end{thm}
\vspace{0.2cm}
In the same article, the author also presented a lower bound on $\dim C_k^r(\Delta)$. The formula can be obtained by replacing \;$\tilde t_i$ by \,$t_i$ in (\ref{upperboundSchumaker}). We shall prove that that formula for the lower bound is the same as the formula we presented in Theorem \ref{lowbHOM}. Essentially the same proof shows that the upper bound formulas, the one in the previous theorem and the one we presented in Theorem (\ref{upbHOM}), coincide. With the exception that the upper bound presented by Schumaker can only be applied for certain numberings on the vertices. That restriction makes sometimes not possible to find an upper bound. In Section \ref{examples} we will include some examples of those situations, and also some cases when not having that restriction leads to find the exact dimension of the space.\\
\\
\textbf{Remark.\;}The above cited article defined $\tilde{t_i}$ as the
number of slopes of the edges containing $\gamma_i$ but not containing
any of the $i-1$ vertices in the list. This coincides with the
definition of $\tilde t_i$ we consider here, except that the reverse
counting on the vertices in $\Delta^0_0$ is used.
\begin{lem}\label{relation-ti}
Let \,$t\geq 2$ and \,$\Omega$ defined as in (\ref{formulas}) above. The notations in the upper bound formulas (and the lower bounds, respectively) are linked in the following way: if\quad$\Omega= r+\ell$\; then\quad$r+j+1-j\cdot t>0$ only when $j\leq\ell-1$.
\end{lem}
\begin{pf}
For any value of $t\geq 2$ we have  \[1<\frac{t}{t-1}\leq 2\,.\] Thus, \;$r+1\leq\Omega\leq 2r+1$, and we can write \;$\Omega=r+\ell$\; for some integer $\ell$\; between 1 and $r+1$. The interval of $t$ where the value of $\Omega$ is $r+\ell$ is given as follows:
\begin{equation}\label{valueOmega}
\begin{cases}
\ST\frac{r+\ell}{\ell}<t\leq\frac{r+(\ell-1)}{\ell-1}\hspace{0.5cm}&\mbox{when\quad$\ell\geq 2$ }\\
\ST\; t>r+1\hspace{0.5cm}&\mbox{when\quad$\ell =1$.}
\end{cases}
\end{equation}
On the other hand, for fixed $r$ and $t$, the condition
\begin{equation}\label{ineq}
r+j+1-j\cdot t\geq 1
\end{equation}
is satisfied if and only if \;$t\leq (r+j)/j$.
Thus, the biggest number $j$ subject to condition (\ref{ineq}) must satisfy
\[\frac{r+(j+1)}{j+1}< t\leq\frac{r+j}{j}\,.\]
From (\ref{valueOmega}), the previous relation holds  if and only if  $\Omega=r+(j+1)$.
\end{pf}
Let \,$t\geq 2$ and \,$\Omega$, $a$ and $b$ be  defined as in (\ref{formulas}). We consider the following two formulas:
{\small
\begin{equation}\label{homverticesformula}
t\binom{k+2-(r+1)}{2}-b\,\binom{k+2-\Omega_i}{2}-a\,\binom{k+2-(\Omega+1)}{2}\,
\end{equation}}\\
{\small
\begin{equation}\label{schverticesformula}
\binom{k+2}{2}-\binom{r+2}{2}-\sum_{j=1}^{k-r}(r+j+1-j\cdot t)_+\,.
\end{equation}}
\begin{lem}\label{ell=1}
If\, $t>r+1$ and $k\geq r$ then the formulas (\ref{homverticesformula}) and (\ref{schverticesformula}) are equal.
\end{lem}
\begin{pf}
By (\ref{valueOmega}) we know that if $t>r+1$ one has $\Omega=r+1$. From the definition in (\ref{formulas}), we have  \,$a=r+1$ and \,$b=t-(r+2)$. Thus, (\ref{homverticesformula}) reduces to
\[{
\begin{cases}
 0&\text{if}\; k+1-r<2\quad(k=r)\\
\ST r+2&\text{if}\; k-r=1\\
 \binom{k+2}{2}-\binom{r+2}{2}&\text{otherwise}.\label{caso1Hom}\\
\end{cases}}
\]
This can be simplified to $\binom{k+2}{2}-\binom{r+2}{2}$ for any value of $k$ and $r$. By Lemma \ref{relation-ti}, the relation (\ref{ineq}) in this case is not satisfied by any value of $j$,  then also (\ref{schverticesformula}) simplifies to that expression.
\end{pf}
\begin{lem}\label{ell>=2}
Let $t$ be an integer $\geq 2$. If \,$t\leq r+1$ and $k\ge r$ then the formulas (\ref{homverticesformula}) and (\ref{schverticesformula}) are equal.
\end{lem}

\begin{pf}
By Lemma \ref{relation-ti} we know that  \,$\Omega=r+\ell$ for some
integer $\ell \geq 2$.  Depending on the value of $\ell$ the binomials in (\ref{homverticesformula}) could be automatically zero. We need to consider three possible situations $k-r\geq\ell+1$, $k-r=1$ and $k-r<\ell$. The formula  reduces to a different expression in each case, but it is not difficult to check that they are respectively equivalent to the expressions we get from (\ref{schverticesformula}).
\end{pf}
\begin{prop}\label{proplower}
The lower bound on $\dim C_k^r(\Delta)$ given in Theorem \ref{lowbHOM}
coincides with the lower bound  in \cite{Sch}[Theorem 3.1 p. 256].
\end{prop}
\begin{pf}
Since any interior vertex of $\Delta$ is contained in at least two edges with different slopes, then $t_i\geq 2$ for every $i=1,2,\dots f_0^0$. Collecting the vertices $\gamma_i$ with the same value $t_i$, the statement follows directly applying Lemmas \ref{ell=1} and \ref{ell>=2} for those values of $t_i$.
\end{pf}
\begin{prop}\label{propupper}
If the vertices in $\Delta^0$ are numbered as in Theorem
(\ref{SchTheorem}), then the upper bound on $\dim C_k^r(\Delta)$ given
in Theorem \ref{upbHOM} coincides with the bound
(\ref{upperboundSchumaker}) in \cite{Sch}[Theorem 2.1 p. 252].
\end{prop}
\begin{pf}
This follows similarly as in the latter proposition. By collecting the vertices in $\Delta^0$ with the same value $\tilde t_i$, we apply Lemmas \ref{ell=1} and \ref{ell>=2} with the values  $\tilde t_i$. The only case that remains to be considered is $\tilde t_i=1$.
It corresponds to the third term in the formula in Theorem (\ref{upbHOM}). Since
{\small
\begin{align*}
\binom{k+2-(r+1)}{2}=\binom{k+2}{2}-\binom{r+2}{2}-\sum_{j=1}^{k-r}(r+1)
\end{align*}}
the statement follows.
\end{pf}

We summarize the results in the following theorem.
Let us denote by  $\textsc {lbh}$, $\textsc {ubh}$ and  $\textsc{lbs}$, $\textsc{ubs}$ the respective lower and upper bounds obtained by Theorems \ref{lowbHOM} and \ref{upbHOM} and the ones obtained in \cite{Sch}.
\begin{thm}Let $\Delta$ be a connected, finite two dimensional simplicial complex, supported on a disk. Then,
\[{\textsc{ lbh}}(\dim C_k^r(\Delta))= \textsc{lbs}(\dim C_k^r(\Delta))\]
\[{\textsc{ubh}}(\dim C_k^r(\Delta))\leq \textsc{ubs}(\dim C_k^r(\Delta))\,.\]
\end{thm}
\begin{pf}
This follows from Propositions \ref{proplower} and \ref{propupper}, and the fact that the formula in Theorem \ref{upbHOM} can be applied to any numbering on the interior vertices of $\Delta$.
\end{pf}
\section{Dimension formula for degree $k\geq 4r+1$}\label{exactdimension}
In this section we present an alternative proof for the dimension
formula of the spline space $C_k^r(\Delta)$ for $k\geq 4r+1$. The
proof is substantially shorter than the one presented in \cite{alfschu}.
We include the following notation: let $\Delta_i$ and
$\Delta_i^{\partial}$ be (respectively) the sets of $i$-dimensional
faces and $i$-dimensional boundary faces. Denote by $f_i(\Delta)$ and
$f_i^{\partial}(\Delta)$ the cardinality of the preceding sets, respectively.
We begin by stating the following lemma which we will need later on.
\begin{lem}\cite[Lemma 3.3]{Local}\label{lema}
If $\Delta$ is a triangulated region in $\R^2$, then there exists a total order on $\Delta_0$ such that for every $\gamma$ in $\Delta^0_0$, there exist vertices $\gamma'$, $\gamma''$ adjacent to $\gamma$, with $\gamma\succ \gamma',\gamma''$, and such that $\overline{\gamma\gamma'}$, $\overline{\gamma\gamma''}$ have distinct slopes.
\end{lem}
For an ordering on $\Delta_0$ as in the previous lemma, we assign indices to the vertices in such a way that $\gamma_i\succ \gamma_j$ when\ $i> j$. The first indices  $1,2,\dots ,f_0^{\partial}$ are assigned to the vertices lying on the boundary. To those interior vertices which are joined to the boundary by two or more edges of distinct slope we assign the indices $f_0^{\partial}+1,\dots,n$.
Let us recall the following notation and remarks from \cite{Local}.
For each interior vertex $\gamma$, and each $f\in\J(\gamma)$, let $f[\gamma]$ denote the corresponding element in $H_0(\J)$. Then $H_0(\J)$ is generated by $\{f[\gamma]\,\mid f\in\J(\gamma)\}$. By definition of $\J(\gamma)$, we know that
\begin{equation}\label{sum}
f[\gamma]=\sum_{\tau\ni \gamma}\ell_\tau^{r+1}f_\tau[\gamma]
\end{equation}
for some polynomials $f_\tau\in R$.
Notice that if $\tau$ is an edge whose vertices are $\gamma$ and $\gamma'$ then the
\begin{equation}\label{opvertexrelation}
\ell_\tau^{r+1}f\,[\gamma]=\ell_\tau^{r+1}[\gamma']
\end{equation}
in $H_0(\J)$; in particular $\ell_\tau^{r+1}f\,[\gamma]=0$ when $\tau$ is an edge connecting $\gamma$ to the boundary.

Here is another lemma that we need to prove that $H_{0} (\J_{k})=0$ for degree $k\geq 4r+1$:
\begin{lem}\label{lem:5.3} Let $\ell_1,\ell_2$ and $\ell_3$ be three equations of
 distinct lines through a point $p$ and $L$ the equation of a line that
 does not contain the point $p$. Then for any polynomial $g$ of degree
 $d \ge r+1$, there exist $u,v \in R$ of degree $d$  and $w\in R$
 of degree $r-1$ such that
$$
\ell_{3}^{r+1} \, g = \ell_{1}^{r+1}\, u + \ell_{2}^{r+1}\, v + L^{d-r+1}\, \ell_{3}^{r+1}\, w,
$$
with $w=0$ in the case $r=0$.
\end{lem}
\begin{pf}
The case $r=0$ is direct from the linear dependency of $\ell_{1}, \ell_{2}$ and $\ell_{3}$.
Suppose that $r>0$. Since $\ell_1$, $\ell_2$ and $L$ are linearly independent, we can make the following change of coordinates:
\begin{align*}
\ell_1&= x\,;&\ell_3&=x+ay\quad\text{for some $a\neq 0$\,;}\\
\ell_2&=y\,;&L&=z\,.
\end{align*}
Let $g=x^i y^j z^k$ be a monomial such that $i+j+k=d$.
If $k\leq d-r$ then $i+j\ge r$ and the polynomial $ \ell_2^{r+1}g$ is
in the ideal $\,\in\langle x^{r+1}, y^{r+1}\rangle$.
If $k\geq d-r+1$ then $g$ is a multiple of $z^{d-r+1}=L^{d-r+1}$. Thus we can write
$$
\ell_{3}^{r+1} \, g = \ell_{1}^{r+1} \, u' + \ell_{2}^{r+1}\, v' + L^{d-r+1} \, w',
$$
for some polynomials $u', v'\in R$ of degree $d$ and \,$w'\in R$ of degree
$2\,r$. As $L^{d-r+1} w'$ is in the ideal $\langle \ell_{1}^{r+1}, \ell_{2}^{r+1},\ell_{3}^{r+1}\rangle$
and as $L$ is a non-zero divisor modulo this ideal, we deduce that
$w'= \ell_{1}^{r+1} \, u'' + \ell_{2}^{r+1}\, v'' + \ell_{3}^{r+1}\, w''$, with
$u'',v'',w'' \in R$ of degree $r-1$.
Collecting the coefficients of $\ell_{1}^{r+1}, \ell_{2}^{r+1}$, we obtain
the desired decomposition:
$$
\ell_{3}^{r+1} \, g = \ell_{1}^{r+1} \, (u' + L^{d-r+1}\,u'') + \ell_{2}^{r+1}\,
(v'+ L^{d-r+1}\, v'') + L^{d-r+1} \,\ell_{3}^{r+1}\, w''.
$$
\end{pf}

We use this lemma to prove the following result:\\
\begin{thm}\label{4r+1}
The dimension of $C_k^{r}(\Delta)$ when $k\geq 4r+1$ is given by the
lower bound formula of Theorem \ref{lowbHOM}.
\end{thm}
\begin{pf}
 From (\ref{eq}), and the formulas for the dimension of the modules
  (\ref{triang}), (\ref{edges}) and (\ref{vertices}), it suffices to
  show that $H_0(\J)_k=0$ for $k\ge 4 r+1$.  Equivalently, we need to show
  that $f[\gamma]=0$ in $H_0(\mathcal{J})_k$ for all $\gamma$ and $f
  \in\J(\gamma)_{k}$. Ordering the vertices as in Lemma \ref{lema}, we
  consider the first interior vertex in the ordering, and denote it by
  $\gamma$. Let $\tau_1, \tau_2, \dots, \tau_t$ be the edges (not
  necessarily with different slopes) containing $\gamma$ and
  $\omega_1$, $\omega_2$, ..., $\omega_t$ be respectively the end
  points of $\tau_{i}$, which are distinct from $\gamma$.  We number
first the edges $\tau_1$ and
  $\tau_2$ that connect $\gamma$ to the boundary, and the remaining
 edges counterclockwise, starting from $\tau_2$.  Let us denote by
  \,$\ell_1,\ell_2,\dots, \ell_t$ (respectively) equations of
  these edges and by $L_{i}$ equations of the edges connecting $\omega_i$ to
$\omega_{i-1}$ for $i=2,\dots t$.

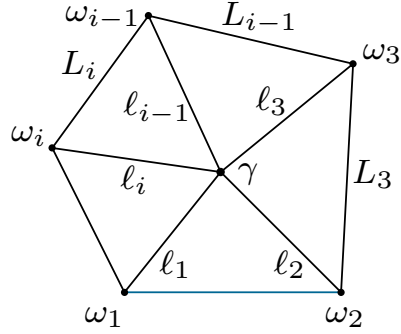
\begin{figure}[!ht]
\scalebox{1.6} 
{\hspace{1.5cm}
\begin{pspicture}(0,-1.549375)(5.490625,1.549375)
\psline[linewidth=0.015cm](2.6734376,-0.13375)(1.8734375,-1.13375)
\psline[linewidth=0.015cm](2.0734375,1.16625)(2.6734376,-0.13375)

\psline[linewidth=0.015cm](1.2734375,0.06625)(2.6734376,-0.13375)

\psline[linewidth=0.015cm](1.8734375,-1.13375)(1.2734375,0.06625)
\psline[linewidth=0.015cm](1.2734375,0.06625)(2.0734375,1.16625)
\psline[linewidth=0.015cm](2.6734376,-0.13375)(3.6734376,-1.13375)
\psline[linewidth=0.015cm,linecolor=MidnightBlue](1.8734375,-1.13375)(3.6734376,-1.13375)
\psline[linewidth=0.015cm](2.6734376,-0.13375)(3.7734375,0.76625)
\psline[linewidth=0.015cm](3.7734375,0.76625)(2.0734375,1.16625)
\psline[linewidth=0.015cm](3.6734376,-1.13375)(3.7734375,0.76625)
\psdots[dotsize=0.06](1.8734375,-1.13375)
\psdots[dotsize=0.06](3.6734376,-1.13375)
\psdots[dotsize=0.06](3.7734375,0.76625)
\psdots[dotsize=0.06](2.0734375,1.16625)
\psdots[dotsize=0.06](1.2734375,0.06625)
\psdots[dotsize=0.07](2.6734376,-0.13375)
\usefont{T1}{ppl}{m}{n}
\rput(2.3,-0.9){\tiny $\ell_1$}
\usefont{T1}{ppl}{m}{n}
\rput(3.25,-0.9){\tiny $\ell_2$}
\usefont{T1}{ppl}{m}{n}
\rput(3.1,0.46125){\tiny $\ell_3$}
\usefont{T1}{ppl}{m}{n}
\rput(2.15,0.4){\tiny $\ell_{i-1}$}
\usefont{T1}{ppl}{m}{n}
\rput(1.95875,-0.2){\tiny $\ell_i$}
\usefont{T1}{ppl}{m}{n}
\rput(3,1.15){\tiny $L_{i-1}$}
\usefont{T1}{ppl}{m}{n}
\rput(1.5,0.76125){\tiny $L_i$}
\usefont{T1}{ppl}{m}{n}
\rput(3.94,-0.15){\tiny $L_3$}
\usefont{T1}{ppl}{m}{n}
\rput(2.9,-0.15){\tiny $\gamma$}
\usefont{T1}{ppl}{m}{n}
\rput(4.02,0.85){\tiny $\omega_3$}
\usefont{T1}{ppl}{m}{n}
\rput(1.7,1.15){\tiny $\omega_{i-1}$}
\usefont{T1}{ppl}{m}{n}
\rput(1.08,0.16125){\tiny $\omega_i$}
\usefont{T1}{ppl}{m}{n}
\rput(1.70875,-1.33875){\tiny $\omega_1$}
\usefont{T1}{ppl}{m}{n}
\rput(3.70875,-1.33875){\tiny $\omega_2$}
\end{pspicture}
}
\caption{Notation in Theorem \ref{4r+1}.}\label{notation near gamma}
\end{figure}

We are going to prove by induction on $j>2$ that for any homogeneous polynomial $f$ of
degree $\ge 3 r$, we have
\begin{align*}
\bullet\;&\ell_{j}^{r+1}\,f\, [\omega_{j}] = \ell_{j}^{r+1} f[\gamma]= 0,\\
\bullet \;& L_{j}^{r+1}\,f [\omega_{j}] =0
\end{align*}
in $H_0(\J)$.

Let us prove it first for $j=3$. Let $f\in R_{k}$ with $k\geq 3r$.
By construction $\omega_2$ is on the boundary, thus we have
$L_{3}^{r+1}f\,[\omega_3]=L_{3}^{r+1}f\,[\omega_2]=0$.
By Lemma \ref{lem:5.3},
$$
\ell_3^{r+1}\,f = \ell_1^{r+1} u + \ell_2^{r+1} v + L_{3}^{2r+1} \ell_{3}^{r+1} w
$$
for some polynomials $u,v, w\in R$.
Then we have
\begin{eqnarray*}
{\ell_3^{r+1}\, f\, [\gamma]}
 & = & \ell_1^{r+1} u\, [\gamma] + \ell_2^{r+1} v \, [\gamma]
+ L_{3}^{2\,r+1} \ell_{3}^{r+1} w \, [\gamma] \\
&=&
L_{3}^{2\,r+1} \ell_{3}^{r+1}\, w\, [\omega_{3}] =  0.
\end{eqnarray*}

Let us take now $i>2$ and assume that the induction hypothesis is true
for $3 \leq j<i$.
We consider first\, $\ell_i^{r+1}f\,[\omega_i]$ with $f$ homogeneous of degree
$\ge 3r$. By Lemma
\ref{lem:5.3} applied to $\ell_{1}, \ell_{2}, \ell_{3}$ and $L_{i}$, we have
$$
\ell_{i}^{r+1}\, f = \ell_{1}^{r+1} \, u + \ell_{2}^{r+1} \, v
+ L_{i}^{2 r+1} \, \ell_{i}^{r+1} w,
$$
for some polynomials $u,v,w\in R$. Then we have
\begin{eqnarray*}
\ell_{i}^{r+1}\, f \, [\gamma] &=& \ell_{1}^{r+1} \, u \, [\gamma]+
\ell_{2}^{r+1} \, v \, [\gamma] + L_{i}^{2 r+1} \, \ell_{i}^{r+1} w \,[\gamma]\\
&=& L_{i}^{2 r+1} \, \ell_{i}^{r+1} w \,[\omega_{i}]
= L_{i}^{2 r+1} \, \ell_{i}^{r+1} w \,[\omega_{i-1}].
\end{eqnarray*}
As $L_i^{2r+1}$ is in the ideal generated by $\ell_{i-1}^{r+1}$, $L_{i-1}^{r+1}$,
by the induction hypothesis, we have
$$
L_{i}^{2 r+1} \, \ell_{i}^{r+1} w \,[\omega_{i-1}]=0.
$$
This proves our claim for $\ell_i^{r+1}\,f\,[\omega_i]=0$.

Let us consider now
$L_i^{r+1}\,f\,[\omega_{i}]=L_i^{r+1}f\,[\omega_{i-1}]$ with $f$
homogeneous of degree $\geq 3r$.
By Lemma \ref{lem:5.3},
$$
L_{i}^{r+1}\, f = \ell_{i-1}^{r+1} \, u + L_{i-1}^{r+1} \, v
+ \ell_{i}^{r+1} \, L_{i}^{r+1}\, w,
$$
for some polynomials $u,v,w\in R$. We deduce that
\begin{eqnarray*}
L_{i}^{r+1}\, f\, [\omega_{i-1}]
&=& \ell_{i-1}^{r+1} \, u\, [\omega_{i-1}]
+ L_{i-1}^{r+1} \, v \, [\omega_{i-1}]
+ \ell_{i}^{r+1} \, L_{i}^{r+1}\, w \, [\omega_{i-1}] \\
&=&
\ell_{i}^{r+1} \, L_{i}^{r+1}\, w \, [\omega_{i}] =
\ell_{i}^{r+1} \, L_{i}^{r+1}\, w \, [\gamma] = 0
\end{eqnarray*}
by the induction hypothesis and the previous computation.

This concludes the induction proof and shows that for any $f \in \J (\gamma)_{4r+1}$,
we have
\begin{align*}
f\, [\gamma]&=\sum_{i=1}^t \ell_i^{r+1}f_i[\gamma]=\sum_{i=3}^t \ell_i^{r+1}f_i[\gamma]\\
& = 0.
\end{align*}
Therefore $H_{0} (\J)=0$ and the dimension of $C^{r}_{k} (\Delta)$ is
given by the lower bound when $k\leq 4r+1$.
\end{pf}

\section{Examples and remarks}\label{examples}
In this section we compute the dimension formula for some complexes. We prove that the dimension formula is in fact given by  the  corresponding lower bound given in Theorem \ref{lowbHOM}. In the first example we show directly that $H_0(\J)_k=0$. In the second example we number the vertices in such a way that the upper bound we obtain  agrees with the lower one. The subdivisions we consider are two of the so--called Powell-Sabin subdivisions (\cite{ps}).

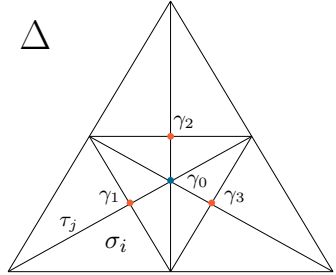
\begin{figure}[ht!]
\centering
\scalebox{0.9} 
{
\begin{pspicture}(0,-2.01)(6.01,2.01)
\psline[linewidth=0.010cm](0.6,-2.0)(3.0,2.0)
\psline[linewidth=0.010cm](3.0,2.0)(5.4,-2.0)
\psline[linewidth=0.010cm](5.4,-2.0)(0.6,-2.0)
\psline[linewidth=0.010cm](0.6,-2.0)(4.2,0.0)
\psline[linewidth=0.010cm](4.2,0.0)(1.8,0.0)
\psline[linewidth=0.010cm](1.8,0.0)(5.4,-2.0)
\psline[linewidth=0.010cm](3.0,2.0)(3.0,-2.0)
\psline[linewidth=0.010cm](3.0,-2.0)(4.2,0.0)
\psline[linewidth=0.010cm](3.0,-2.0)(1.8,0.0)
\usefont{T1}{ppl}{m}{n}
\rput(3.4,-0.67){\scriptsize $\gamma_0$}
\usefont{T1}{ppl}{m}{n}
\rput(2.1,-0.9){\scriptsize $\gamma_1$}
\usefont{T1}{ppl}{m}{n}
\rput(3.2,0.2){\scriptsize $\gamma_2$}
\usefont{T1}{ppl}{m}{n}
\rput(3.95,-0.9){\scriptsize $\gamma_3$}
\usefont{T1}{ppl}{m}{n}
\rput(1,1.5){\Large $\Delta$}
\usefont{T1}{ppl}{m}{n}
\rput(1.5,-1.3){\scriptsize $\tau_j$}
\usefont{T1}{ppl}{m}{n}
\rput(2.2,-1.6){\small $\sigma_i$}
\psdots[dotsize=0.1,linecolor=MidnightBlue](3.0,-0.65)
\psdots[dotsize=0.1,linecolor=RedOrange](2.39999,-0.984)
\psdots[dotsize=0.1,linecolor=RedOrange](3.609,-0.984)
\psdots[dotsize=0.1,linecolor=RedOrange](3.0,0.0)
\end{pspicture}
}
\caption{Powell-Sabin 12-split.}\label{12split}
\end{figure}

\begin{example} Powell-Sabin 12-split.
\end{example}
Let   $\Delta$  be the simplicial complex supported on a triangle $|\Delta|$ subdivided into twelve smaller triangles, as in Figure  \ref{12split}.

Let us note that by using any numbering on $\Delta_0^0$, the upper bound we get from Theorem \ref{upbHOM} equals the lower bound, leading easily to the dimension formula. On the other hand, it is not possible to find a numbering on the interior vertices of $\Delta$ with the condition required to apply Theorem \ref{SchTheorem}. Hence it is not possible to find an upper bound for this spline space by using the formula of Schumaker.

By the remarks in Section \ref{exactdimension}, it is easy to check that $ \ell_j^{r+1}f\,[\gamma_i]=0$ for every $\gamma_i$ and $\ell_j^{r+1}f\,\in \J(\gamma_i)$. Thus, for any $k$ and any $r\leq k$, the homology module $H_0(\J)_k$ is zero and the dimension of the spline space $C_k^r(\Delta)$ is given by the formula in Theorem \ref{lowbHOM}.

This example is an instance of a triangulation where all
edges are pseudoboundaries (in the terminology of \cite{family}), and the
computation also follows by Lemma 2.5 of \cite{family}.

\begin{example} Powell-Sabin 6-split.
\end{example}
Let $\Delta$ be the simplicial complex supported on a simply connected triangulated region $|\Delta|$ in $\R^2$. Assume the triangulation consists of $f_0$ vertices (interior and exterior) and $f_i^0$ $i$-dimensional interior faces ($i=0,1,2$), where $f_2^0$ is taken as the number of all the triangles $\sigma_i$ in the subdivision.
We refine this triangulation  by subdividing each triangle into 6 triangles to get $\tilde\Delta$, as follows.
\begin{enumerate}
\item In each triangle $\sigma_i\in\Delta^0_2$ choose an interior point $\nu_i$, in such a way that if two triangles \,$\sigma_i$, $\sigma_j$ have an edge $\tau$ in common, then the line joining $\nu_i$  and $\nu_j$ intersects $\tau$ at an interior point $\mu_{ij}$.
\item Join each new point $\nu_i$ to the vertices of the triangle $\sigma_i$, and to the points $\mu_{ij}$ on the edges of $\sigma_i$.
\item For the triangles having an edge (or more) on the boundary, choose a point in the interior of such edge and join it with $\nu_i$.
\end{enumerate}

\begin{figure}[!ht]
\scalebox{0.8} 
{\begin{pspicture}(-2,10.063)(26.71,6.399)
\definecolor{RoyalPurple}{rgb}{0.28627450980392155,0.21176470588235294,0.3411764705882353}
\definecolor{color2746}{rgb}{0.06666666666666667,0.4627450980392157,0.2}
\definecolor{color2803}{rgb}{0.3843137254901961,0.6901960784313725,0.03137254901960784}
\psdots[dotsize=0.11,linecolor=gray](4.641875,8.7)
\psdots[dotsize=0.11,linecolor=gray](2.741875,7.8)
\psdots[dotsize=0.11,linecolor=gray](3.141875,6.6)
\psdots[dotsize=0.11,linecolor=gray](1.541875,5.9)
\psdots[dotsize=0.11,linecolor=gray](1.941875,4.8)
\psdots[dotsize=0.11,linecolor=gray](3.241875,4.4)
\psdots[dotsize=0.11,linecolor=gray](4.841875,5.0)
\psdots[dotsize=0.11,linecolor=gray](6.341875,4.2)
\psdots[dotsize=0.11,linecolor=gray](7.241875,5.5)
\psdots[dotsize=0.11,linecolor=gray](8.641875,6.2)
\psdots[dotsize=0.11,linecolor=gray](9.941875,5.8)
\psdots[dotsize=0.11,linecolor=gray](10.341875,4.5)
\psdots[dotsize=0.11,linecolor=gray](10.541875,7.4)
\psdots[dotsize=0.11,linecolor=gray](10.141875,8.3)
\psdots[dotsize=0.11,linecolor=gray](8.941875,8.2)
\psdots[dotsize=0.11,linecolor=gray](8.041875,9.0)
\psdots[dotsize=0.11,linecolor=gray](6.441875,8.0)
\psdots[dotsize=0.11,linecolor=gray](5.941875,6.9)
\psdots[dotsize=0.11,linecolor=gray](4.741875,6.4)
\psline[linewidth=0.010cm,linecolor=gray](27.741875,-9.6)(27.641874,-9.4)
\psline[linewidth=0.010cm,linecolor=gray](1.541875,5.9)(1.941875,4.8)
\psline[linewidth=0.010cm,linecolor=gray](2.741875,7.8)(4.641875,8.7)
\psline[linewidth=0.010cm,linecolor=gray](4.641875,8.7)(6.441875,8.0)
\psline[linewidth=0.010cm,linecolor=gray](6.441875,8.0)(8.041875,9.0)
\psline[linewidth=0.010cm,linecolor=gray](8.041875,9.0)(8.941875,8.2)
\psline[linewidth=0.010cm,linecolor=gray](8.941875,8.2)(10.141875,8.3)
\psline[linewidth=0.010cm,linecolor=gray](10.141875,8.3)(10.541875,7.4)
\psline[linewidth=0.010cm,linecolor=gray](10.541875,7.4)(9.941875,5.8)
\psline[linewidth=0.010cm,linecolor=gray](9.941875,5.8)(10.341875,4.5)
\psline[linewidth=0.010cm,linecolor=gray](1.941875,4.8)(3.241875,4.4)
\psline[linewidth=0.010cm,linecolor=gray](3.241875,4.4)(4.841875,5.0)
\psline[linewidth=0.010cm,linecolor=gray](4.841875,5.0)(6.341875,4.2)
\psline[linewidth=0.010cm,linecolor=gray](6.341875,4.2)(7.241875,5.5)
\psline[linewidth=0.010cm,linecolor=gray](7.241875,5.5)(8.641875,6.2)
\psline[linewidth=0.010cm,linecolor=gray](8.641875,6.2)(9.941875,5.8)
\psline[linewidth=0.010cm,linecolor=gray](2.741875,7.8)(3.141875,6.6)
\psline[linewidth=0.010cm,linecolor=gray](3.141875,6.6)(4.741875,6.4)
\psline[linewidth=0.010cm,linecolor=gray](4.741875,6.4)(5.941875,6.9)
\psline[linewidth=0.010cm,linecolor=gray](5.941875,6.9)(6.441875,8.0)
\psline[linewidth=0.010cm,linecolor=gray](5.941875,6.9)(7.241875,5.5)
\psline[linewidth=0.010cm,linecolor=gray](4.741875,6.4)(4.841875,5.0)
\psline[linewidth=0.010cm,linecolor=gray](1.541875,5.9)(3.141875,6.6)
\psline[linewidth=0.010cm,linecolor=gray](8.641875,6.2)(8.941875,8.2)
\psline[linewidth=0.010cm,linecolor=gray](5.941875,6.9)(7.441875,7.3)
\psline[linewidth=0.010cm,linecolor=gray](5.941875,6.9)(4.741875,7.5)
\psline[linewidth=0.010cm,linecolor=gray](5.941875,6.9)(6.041875,5.7)
\psline[linewidth=0.010cm,linecolor=gray](1.141875,6.6)(2.741875,7.8)
\psline[linewidth=0.010cm,linecolor=gray](2.741875,7.8)(2.541875,9.1)
\psline[linewidth=0.010cm,linecolor=gray](2.741875,7.8)(4.741875,7.5)
\psline[linewidth=0.010cm,linecolor=gray](4.741875,7.5)(4.641875,8.7)
\psline[linewidth=0.010cm,linecolor=gray](4.641875,8.7)(7.241875,9.6)
\psline[linewidth=0.010cm,linecolor=gray](4.641875,8.7)(2.541875,9.1)
\psline[linewidth=0.010cm,linecolor=gray](4.741875,7.5)(6.441875,8.0)
\psline[linewidth=0.010cm,linecolor=gray](6.441875,8.0)(7.241875,9.6)
\psline[linewidth=0.010cm,linecolor=gray](6.441875,8.0)(7.441875,7.3)
\psline[linewidth=0.010cm,linecolor=gray](7.441875,7.3)(8.041875,9.0)
\psline[linewidth=0.010cm,linecolor=gray](8.041875,9.0)(9.841875,9.5)
\psline[linewidth=0.010cm,linecolor=gray](7.241875,9.6)(8.041875,9.0)
\psline[linewidth=0.010cm,linecolor=gray](7.441875,7.3)(8.941875,8.2)
\psline[linewidth=0.010cm,linecolor=gray](8.941875,8.2)(9.841875,7.3)
\psline[linewidth=0.010cm,linecolor=gray](9.841875,7.3)(10.141875,8.3)
\psline[linewidth=0.010cm,linecolor=gray](10.141875,8.3)(10.641875,8.2)
\psline[linewidth=0.010cm,linecolor=gray](10.641875,8.2)(10.541875,7.4)
\psline[linewidth=0.010cm,linecolor=gray](10.541875,7.4)(9.841875,7.3)
\psline[linewidth=0.010cm,linecolor=gray](10.541875,7.4)(11.141875,6.6)
\psline[linewidth=0.010cm,linecolor=gray](11.141875,6.6)(9.941875,5.8)
\psline[linewidth=0.010cm,linecolor=gray](9.941875,5.8)(9.841875,7.3)
\psline[linewidth=0.010cm,linecolor=gray](9.841875,7.3)(8.641875,6.2)
\psline[linewidth=0.010cm,linecolor=gray](8.941875,8.2)(9.841875,9.5)
\psline[linewidth=0.010cm,linecolor=gray](9.841875,9.5)(10.141875,8.3)
\psline[linewidth=0.010cm,linecolor=gray](7.441875,7.3)(8.641875,6.2)
\psline[linewidth=0.010cm,linecolor=gray](7.241875,5.5)(7.441875,7.3)
\psline[linewidth=0.010cm,linecolor=gray](6.041875,5.7)(7.241875,5.5)
\psline[linewidth=0.010cm,linecolor=gray](1.141875,6.6)(3.141875,6.6)
\psline[linewidth=0.010cm,linecolor=gray](3.141875,6.6)(4.741875,7.5)
\psline[linewidth=0.010cm,linecolor=gray](3.141875,6.6)(3.341875,5.8)
\psline[linewidth=0.010cm,linecolor=gray](3.341875,5.8)(4.741875,6.4)
\psline[linewidth=0.010cm,linecolor=gray](4.741875,6.4)(4.741875,7.5)
\psline[linewidth=0.010cm,linecolor=gray](4.741875,6.4)(6.041875,5.7)
\psline[linewidth=0.010cm,linecolor=gray](6.041875,5.7)(4.841875,5.0)
\psline[linewidth=0.010cm,linecolor=gray](4.841875,5.0)(3.341875,5.8)
\psline[linewidth=0.010cm,linecolor=gray](4.841875,5.0)(5.241875,3.5)
\psline[linewidth=0.010cm,linecolor=gray](5.241875,3.5)(6.341875,4.2)
\psline[linewidth=0.010cm,linecolor=gray](6.041875,5.7)(6.341875,4.2)
\psline[linewidth=0.010cm,linecolor=gray](6.341875,4.2)(8.941875,3.5)
\psline[linewidth=0.010cm,linecolor=gray](8.641875,6.2)(8.941875,3.5)
\psline[linewidth=0.010cm,linecolor=gray](9.941875,5.8)(8.941875,3.5)
\psline[linewidth=0.010cm,linecolor=gray](10.341875,4.5)(8.941875,3.5)
\psline[linewidth=0.010cm,linecolor=gray](10.341875,4.5)(10.641875,4.1)
\psline[linewidth=0.010cm,linecolor=gray](10.341875,4.5)(11.141875,6.6)
\psline[linewidth=0.010cm,linecolor=gray](7.241875,5.5)(8.941875,3.5)
\psline[linewidth=0.010cm,linecolor=gray](3.241875,4.4)(1.741875,3.5)
\psline[linewidth=0.010cm,linecolor=gray](3.241875,4.4)(5.241875,3.5)
\psline[linewidth=0.010cm,linecolor=gray](3.241875,4.4)(3.341875,5.8)
\psline[linewidth=0.010cm,linecolor=gray](3.341875,5.8)(1.941875,4.8)
\psline[linewidth=0.010cm,linecolor=gray](1.941875,4.8)(1.741875,3.5)
\psline[linewidth=0.010cm,linecolor=gray](1.941875,4.8)(1.041875,4.9)
\psline[linewidth=0.010cm,linecolor=gray](1.041875,4.9)(1.541875,5.9)
\psline[linewidth=0.010cm,linecolor=gray](1.541875,5.9)(3.341875,5.8)
\psline[linewidth=0.010cm,linecolor=gray](1.541875,5.9)(1.141875,6.6)
\psline[linewidth=0.010cm,linecolor=gray](4.641875,8.7)(4.541875,9.3)
\psline[linewidth=0.010cm,linecolor=gray](2.741875,7.8)(1.941875,8.0)
\psline[linewidth=0.010cm,linecolor=gray](8.041875,9.0)(8.041875,9.6)
\psline[linewidth=0.010cm,linecolor=gray](10.141875,8.3)(10.441875,8.5)
\psline[linewidth=0.010cm,linecolor=gray](10.541875,7.4)(10.91,7.4)
\psline[linewidth=0.010cm,linecolor=gray](10.341875,4.5)(10.741875,4.5)
\psline[linewidth=0.010cm,linecolor=gray](10.341875,4.5)(10.141875,3.9)
\psline[linewidth=0.010cm,linecolor=gray](3.241875,4.4)(3.141875,3.5)
\psline[linewidth=0.010cm,linecolor=gray](1.541875,5.9)(1.09,5.8)
\psline[linewidth=0.010cm,linecolor=gray](1.341875,4.3)(1.941875,4.8)
\psline[linewidth=0.010cm,linecolor=gray](6.341875,4.2)(6.441875,3.5)
\psline[linewidth=0.040cm](8.941875,3.5)(11.141875,6.6)
\psline[linewidth=0.040cm](7.441875,7.3)(8.941875,3.5)
\psline[linewidth=0.040cm](8.941875,3.5)(9.841875,7.3)
\psline[linewidth=0.040cm](6.041875,5.7)(8.941875,3.5)
\psline[linewidth=0.040cm](6.041875,5.7)(5.241875,3.5)
\psline[linewidth=0.040cm](5.241875,3.5)(3.341875,5.8)
\psline[linewidth=0.040cm](1.741875,3.5)(3.341875,5.8)
\psline[linewidth=0.040cm](9.841875,7.3)(11.141875,6.6)
\psline[linewidth=0.040cm](7.441875,7.3)(9.841875,7.3)
\psline[linewidth=0.040cm](7.441875,7.3)(6.041875,5.7)
\psline[linewidth=0.040cm](3.341875,5.8)(6.041875,5.7)
\psline[linewidth=0.040cm](1.741875,3.5)(1.041875,4.9)
\psline[linewidth=0.040cm](1.041875,4.9)(3.341875,5.8)
\psline[linewidth=0.040cm](1.041875,4.9)(1.141875,6.6)
\psline[linewidth=0.040cm](1.141875,6.6)(3.341875,5.8)
\psline[linewidth=0.040cm](3.341875,5.8)(4.741875,7.5)
\psline[linewidth=0.040cm](4.741875,7.5)(1.141875,6.6)
\psline[linewidth=0.040cm](1.141875,6.6)(2.541875,9.1)
\psline[linewidth=0.040cm](2.541875,9.1)(4.741875,7.5)
\psline[linewidth=0.040cm](4.741875,7.5)(6.041875,5.7)
\psline[linewidth=0.040cm](4.741875,7.5)(7.441875,7.3)
\psline[linewidth=0.040cm](7.441875,7.3)(9.841875,9.5)
\psline[linewidth=0.040cm](9.841875,9.5)(9.841875,7.3)
\psline[linewidth=0.040cm](9.841875,7.3)(10.641875,8.2)
\psline[linewidth=0.040cm](11.141875,6.6)(10.641875,8.2)
\psline[linewidth=0.040cm](10.641875,8.2)(9.841875,9.5)
\psline[linewidth=0.040cm](7.441875,7.3)(7.241875,9.6)
\psline[linewidth=0.040cm](7.241875,9.6)(9.841875,9.5)
\psline[linewidth=0.040cm](4.741875,7.5)(7.241875,9.6)
\psline[linewidth=0.040cm](2.541875,9.1)(7.241875,9.6)
\psline[linewidth=0.040cm](1.741875,3.5)(5.241875,3.5)
\psline[linewidth=0.040cm](5.241875,3.5)(8.941875,3.5)
\psline[linewidth=0.040cm](8.941875,3.5)(10.641875,4.1)
\psline[linewidth=0.040cm](10.641875,4.1)(11.141875,6.6)
\usefont{T1}{ppl}{m}{n}
\rput(1.35,9.3){\Huge{$\tilde\Delta$}}
\end{pspicture}
}
\vspace{2cm}
\caption{Powell-Sabin 6-split of $\Delta$.}\label{6splitcomplete}
\end{figure}
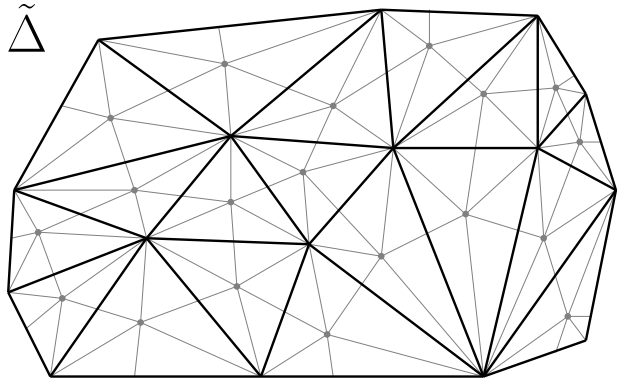

We want to apply Corollary \ref{>r+1} and give a formula for $\dim C_k^1(\tilde\Delta)$.

\begin{figure}[!ht]
\scalebox{0.8} 
{\hspace{0.6cm}
\begin{pspicture}(0,-9.659375)(29.07375,9.694375)
\definecolor{RoyalPurple}{rgb}{0.28627450980392155,0.21176470588235294,0.3411764705882353}
\definecolor{color142}{rgb}{0.9098039215686274,0.1803921568627451,0.06666666666666667}
\definecolor{color82}{rgb}{0.06666666666666667,0.4627450980392157,0.2}
\psline[linewidth=0.03cm,linecolor=gray](8.76875,7.245625)(10.26875,3.445625)
\psline[linewidth=0.03cm,linecolor=gray](10.26875,3.445625)(11.16875,7.245625)
\psline[linewidth=0.03cm,linecolor=gray](8.76875,7.245625)(11.16875,7.245625)
\psline[linewidth=0.03cm,linecolor=gray](8.76875,7.245625)(7.36875,5.645625)
\psline[linewidth=0.03cm,linecolor=gray](4.66875,5.745625)(7.36875,5.645625)
\psline[linewidth=0.03cm,linecolor=gray](4.66875,5.745625)(6.06875,7.445625)
\psline[linewidth=0.03cm,linecolor=gray](6.06875,7.445625)(7.36875,5.645625)
\psline[linewidth=0.03cm,linecolor=gray](6.06875,7.445625)(8.76875,7.245625)
\psline[linewidth=0.02cm,linecolor=RoyalPurple](29.06875,-9.654375)(28.968748,-9.454375)
\psline[linewidth=0.02cm,linecolor=RoyalPurple](2.86875,5.845625)(3.26875,4.745625)
\psline[linewidth=0.02cm,linecolor=RoyalPurple](4.06875,7.745625)(5.96875,8.645625)
\psline[linewidth=0.02cm,linecolor=RoyalPurple](5.96875,8.645625)(7.03,8.25)
\psline[linewidth=0.01cm,linecolor=gray](7.03,8.25)(7.76875,7.945625)
\psdots[dotsize=0.1,linecolor=RoyalPurple](7.03,8.25)%
\psline[linewidth=0.01cm,linecolor=gray](7.76875,7.945625)(8.66,8.494375)
\psline[linewidth=0.02cm,linecolor=RoyalPurple](8.66,8.494375)(9.36875,8.945625)
\psline[linewidth=0.02cm,linecolor=RoyalPurple](9.36875,8.945625)(10.015,8.394375)
\psline[linewidth=0.01cm,linecolor=gray](10.015,8.394375)(10.26875,8.145625)
\psline[linewidth=0.01cm,linecolor=gray](10.26875,8.145625)(11.18,8.23)
\psline[linewidth=0.02cm,linecolor=RoyalPurple](11.18,8.23)(11.46875,8.245625)
\psline[linewidth=0.02cm,linecolor=RoyalPurple](11.46875,8.245625)(11.86875,7.345625)
\psline[linewidth=0.01cm,linecolor=gray](11.73,6.95)(11.26875,5.745625)
\psline[linewidth=0.02cm,linecolor=RoyalPurple](11.73,6.95)(11.86875,7.345625)
\psline[linewidth=0.01cm,linecolor=gray](11.26875,5.745625)(11.47,5.125)
\psline[linewidth=0.02cm,linecolor=RoyalPurple](11.47,5.125)(11.66875,4.445625)
\psline[linewidth=0.02cm,linecolor=RoyalPurple](3.26875,4.745625)(4.56875,4.345625)
\psline[linewidth=0.02cm,linecolor=RoyalPurple](4.56875,4.345625)(5.535,4.699)
\psline[linewidth=0.01cm,linecolor=gray](5.535,4.699)(6.16875,4.945625)
\psline[linewidth=0.01cm,linecolor=gray](6.16875,4.945625)(6.96,4.53)
\psline[linewidth=0.02cm,linecolor=RoyalPurple](6.96,4.53)(7.66875,4.145625)
\psline[linewidth=0.02cm,linecolor=RoyalPurple](7.66875,4.145625)(8.26,4.98)
\psline[linewidth=0.01cm,linecolor=gray](8.26,4.98)(8.56875,5.445625)
\psline[linewidth=0.01cm,linecolor=gray](8.56875,5.445625)(9.96875,6.145625)
\psline[linewidth=0.01cm,linecolor=gray](9.96875,6.145625)(11.26875,5.745625)
\psline[linewidth=0.01cm,linecolor=gray](4.46875,6.545625)(6.06875,6.345625)
\psline[linewidth=0.01cm,linecolor=gray](6.06875,6.345625)(7.26875,6.845625)
\psline[linewidth=0.01cm,linecolor=gray](7.26875,6.845625)(7.76875,7.945625)
\psline[linewidth=0.01cm,linecolor=gray](7.26875,6.845625)(8.56875,5.445625)
\psline[linewidth=0.01cm,linecolor=gray](6.06875,6.345625)(6.16875,4.945625)
\psline[linewidth=0.01cm,linecolor=gray](3.56,6.15)(4.46875,6.545625)
\psline[linewidth=0.01cm,linecolor=gray](9.96875,6.145625)(10.26875,8.145625)
\psline[linewidth=0.01cm,linecolor=gray](7.26875,6.845625)(8.76875,7.245625)
\psline[linewidth=0.01cm,linecolor=gray](7.26875,6.845625)(6.06875,7.445625)
\psline[linewidth=0.01cm,linecolor=gray](7.26875,6.845625)(7.36875,5.645625)
\psline[linewidth=0.01cm,linecolor=gray](6.06875,7.445625)(7.76875,7.945625)
\psline[linewidth=0.01cm,linecolor=gray](4.46875,6.545625)(6.06875,7.445625)
\psline[linewidth=0.01cm,linecolor=gray](6.06875,6.345625)(6.06875,7.445625)
\psline[linewidth=0.01cm,linecolor=gray](8.56875,5.445625)(8.76875,7.245625)
\psline[linewidth=0.01cm,linecolor=gray](7.36875,5.645625)(8.56875,5.445625)
\psline[linewidth=0.01cm,linecolor=gray](8.56875,5.445625)(10.26875,3.445625)
\psline[linewidth=0.01cm,linecolor=gray](7.76875,7.945625)(8.76875,7.245625)
\psline[linewidth=0.01cm,linecolor=gray](8.76875,7.245625)(10.26875,8.145625)
\psline[linewidth=0.01cm,linecolor=gray](8.76875,7.245625)(9.96875,6.145625)
\psline[linewidth=0.01cm,linecolor=gray](7.76875,7.945625)(8.56875,9.545625)
\psline[linewidth=0.01cm,linecolor=gray](10.26875,8.145625)(11.16875,7.245625)
\psline[linewidth=0.01cm,linecolor=gray](10.26875,8.145625)(11.16875,9.445625)
\psline[linewidth=0.01cm,linecolor=gray](11.26875,5.745625)(11.16875,7.245625)
\psline[linewidth=0.01cm,linecolor=gray](11.16875,7.245625)(9.96875,6.145625)
\psline[linewidth=0.01cm,linecolor=gray](12.46875,6.545625)(11.26875,5.745625)
\psline[linewidth=0.01cm,linecolor=gray](11.26875,5.745625)(10.26875,3.445625)
\psline[linewidth=0.01cm,linecolor=gray](7.36875,5.645625)(6.16875,4.945625)
\psline[linewidth=0.01cm,linecolor=gray](6.16875,4.945625)(4.66875,5.745625)
\psline[linewidth=0.01cm,linecolor=gray](6.16875,4.945625)(6.56875,3.445625)
\psline[linewidth=0.01cm,linecolor=gray](4.66875,5.745625)(6.06875,6.345625)
\psline[linewidth=0.01cm,linecolor=gray](4.46875,6.545625)(4.66875,5.745625)
\psline[linewidth=0.01cm,linecolor=gray](9.96875,6.145625)(10.26875,3.445625)
\psline[linewidth=0.01cm,linecolor=gray](4.32,7)(4.46875,6.545625)
\psline[linewidth=0.02cm,linecolor=RoyalPurple](4.06875,7.745625)(4.32,7)
\psline[linewidth=0.02cm,linecolor=MidnightBlue](2.46875,6.545625)(4.06875,7.745625)
\psline[linewidth=0.02cm,linecolor=MidnightBlue](4.06875,7.745625)(3.86875,9.045625)
\psline[linewidth=0.02cm,linecolor=color82](4.06875,7.745625)(6.06875,7.445625)
\psline[linewidth=0.02cm,linecolor=color82](6.06875,7.445625)(5.96875,8.645625)
\psline[linewidth=0.02cm,linecolor=MidnightBlue](8.56875,9.545625)(9.36875,8.945625)
\psline[linewidth=0.04cm,linecolor=color82](6.06875,7.445625)(5.03,8.2)
\psline[linewidth=0.04cm,linecolor=RoyalPurple](3.86875,9.045625)(5.03,8.2)
\psline[linewidth=0.04cm,linecolor=color82](6.06875,7.445625)(4.32,7)
\psline[linewidth=0.04cm,linecolor=RoyalPurple](4.32,7)(2.46875,6.545625)
\psline[linewidth=0.04cm,linecolor=color82](6.06875,7.445625)(7.03,8.25)
\psline[linewidth=0.04cm,linecolor=RoyalPurple](7.03,8.25)(8.56875,9.545625)
\psline[linewidth=0.04cm,linecolor=RoyalPurple](8.66,8.494375)(8.56875,9.545625)
\psline[linewidth=0.04cm,linecolor=color82](8.76875,7.245625)(8.66,8.494375)
\psline[linewidth=0.04cm,linecolor=color82](8.76875,7.245625)(10.015,8.394375)
\psline[linewidth=0.04cm,linecolor=RoyalPurple](10.015,8.394375)(11.16875,9.445625)
\psline[linewidth=0.04cm,linecolor=color82](11.175,8.23)(11.16875,7.245625)
\psline[linewidth=0.04cm,linecolor=RoyalPurple](11.16875,9.445625)(11.175,8.23)
\psline[linewidth=0.04cm,linecolor=RoyalPurple](11.73,6.95)(12.46875,6.545625)
\psline[linewidth=0.04cm,linecolor=color82](11.16875,7.245625)(11.73,6.95)
\psline[linewidth=0.04cm,linecolor=RoyalPurple](11.6663,7.81)(11.96875,8.145625)
\psline[linewidth=0.04cm,linecolor=color82](11.16875,7.245625)(11.6663,7.81)
\psline[linewidth=0.04cm,linecolor=RoyalPurple](10.26875,3.445625)(12.46875,6.545625)
\psline[linewidth=0.04cm,linecolor=RoyalPurple](8.26,4.98)(10.26875,3.445625)
\psline[linewidth=0.04cm,linecolor=color82](7.36875,5.645625)(8.26,4.98)
\psline[linewidth=0.04cm,linecolor=RoyalPurple](6.96,4.53)(6.56875,3.445625)
\psline[linewidth=0.04cm,linecolor=color82](7.36875,5.645625)(6.96,4.53)
\psline[linewidth=0.02cm,linecolor=MidnightBlue](5.96875,8.645625)(8.56875,9.545625)
\psline[linewidth=0.02cm,linecolor=MidnightBlue](5.96875,8.645625)(3.86875,9.045625)
\psline[linewidth=0.02cm,linecolor=MidnightBlue](9.36875,8.945625)(11.16875,9.445625)
\psline[linewidth=0.02cm,linecolor=MidnightBlue](11.16875,9.445625)(11.46875,8.245625)
\psline[linewidth=0.02cm,linecolor=color82](11.16875,7.245625)(11.46875,8.245625)
\psline[linewidth=0.02cm,linecolor=color82](7.36875,5.645625)(7.66875,4.145625)
\psline[linewidth=0.02cm,linecolor=MidnightBlue](11.46875,8.245625)(11.96875,8.145625)
\psline[linewidth=0.02cm,linecolor=MidnightBlue](11.96875,8.145625)(11.86875,7.345625)
\psline[linewidth=0.02cm,linecolor=color82](11.86875,7.345625)(11.16875,7.245625)
\psline[linewidth=0.02cm,linecolor=MidnightBlue](11.86875,7.345625)(12.46875,6.545625)
\psline[linewidth=0.01cm,linecolor=gray](2.46875,6.545625)(4.46875,6.545625)
\psline[linewidth=0.01cm,linecolor=gray](6.06875,6.345625)(7.36875,5.645625)
\psline[linewidth=0.02cm,linecolor=MidnightBlue](6.56875,3.445625)(7.66875,4.145625)
\psline[linewidth=0.02cm,linecolor=MidnightBlue](7.66875,4.145625)(10.26875,3.445625)
\psline[linewidth=0.02cm,linecolor=MidnightBlue](11.66875,4.445625)(10.26875,3.445625)
\psline[linewidth=0.02cm,linecolor=MidnightBlue](11.66875,4.445625)(11.96875,4.045625)
\psline[linewidth=0.02cm,linecolor=MidnightBlue](11.66875,4.445625)(12.46875,6.545625)
\psline[linewidth=0.02cm,linecolor=MidnightBlue](4.56875,4.345625)(3.06875,3.445625)
\psline[linewidth=0.02cm,linecolor=MidnightBlue](4.56875,4.345625)(6.56875,3.445625)
\psline[linewidth=0.02cm,linecolor=color82](4.56875,4.345625)(4.66875,5.745625)
\psline[linewidth=0.02cm,linecolor=color82](4.66875,5.745625)(3.26875,4.745625)
\psline[linewidth=0.02cm,linecolor=MidnightBlue](3.26875,4.745625)(3.06875,3.445625)
\psline[linewidth=0.02cm,linecolor=MidnightBlue](3.26875,4.745625)(2.36875,4.845625)
\psline[linewidth=0.02cm,linecolor=MidnightBlue](2.36875,4.845625)(2.86875,5.845625)
\psline[linewidth=0.02cm,linecolor=color82](2.86875,5.845625)(4.66875,5.745625)
\psline[linewidth=0.02cm,linecolor=MidnightBlue](2.86875,5.845625)(2.46875,6.545625)
\psline[linewidth=0.02cm,linecolor=MidnightBlue](5.96875,8.645625)(5.86875,9.245625)
\psline[linewidth=0.02cm,linecolor=MidnightBlue](4.06875,7.745625)(3.26875,7.945625)
\psline[linewidth=0.02cm,linecolor=MidnightBlue](9.36875,8.945625)(9.36875,9.545625)
\psline[linewidth=0.02cm,linecolor=MidnightBlue](11.46875,8.245625)(11.76875,8.445625)
\psline[linewidth=0.02cm,linecolor=MidnightBlue](11.86875,7.345625)(12.24,7.345625)
\psline[linewidth=0.02cm,linecolor=MidnightBlue](11.66875,4.445625)(12.06875,4.445625)
\psline[linewidth=0.02cm,linecolor=MidnightBlue](11.66875,4.445625)(11.46875,3.845625)
\psline[linewidth=0.02cm,linecolor=MidnightBlue](4.56875,4.345625)(4.46875,3.445625)
\psline[linewidth=0.02cm,linecolor=MidnightBlue](2.86875,5.845625)(2.44,5.794375)
\psline[linewidth=0.02cm,linecolor=MidnightBlue](2.66875,4.245625)(3.26875,4.745625)
\psline[linewidth=0.02cm,linecolor=MidnightBlue](7.66875,4.145625)(7.76875,3.445625)
\psline[linewidth=0.02cm,linecolor=color82](8.76875,7.245625)(9.36875,8.945625)
\psline[linewidth=0.04cm,linecolor=color82](4.66875,5.745625)(3.56,6.15)
\psline[linewidth=0.04cm,linecolor=color82](4.66875,5.745625)(3.12,5.14)
\psline[linewidth=0.04cm,linecolor=color82](4.66875,5.745625)(3.85,4.56)
\psline[linewidth=0.04cm,linecolor=color82](4.66875,5.745625)(5.535,4.699)
\psdots[dotsize=0.1,linecolor=color82](4.66875,5.745625)
\usefont{T1}{ppl}{m}{n}
\rput(7.65,5.2){\small ${ 4}$}
\usefont{T1}{ppl}{m}{n}
\rput(7.2,4.7){\small ${ 2}$}
\usefont{T1}{ppl}{m}{n}
\rput(8.3,4.7){\small ${ 3}$}
\usefont{T1}{ppl}{m}{n}
\rput(7.85,3.9){\small ${1 }$}
\psline[linewidth=0.04cm,linecolor=RoyalPurple](6.56875,3.445625)(5.535,4.699)
\psline[linewidth=0.04cm,linecolor=RoyalPurple](3.06875,3.445625)(3.85,4.56)
\psline[linewidth=0.04cm,linecolor=RoyalPurple](2.46875,6.545625)(3.56,6.15)
\psline[linewidth=0.04cm,linecolor=RoyalPurple](2.36875,4.845625)(3.12,5.14)
\psdots[dotsize=0.1,linecolor=RoyalPurple](4.32,7)
\psdots[dotsize=0.1,linecolor=RoyalPurple](3.56,6.15)
\psdots[dotsize=0.1,linecolor=RoyalPurple](3.12,5.14)
\psdots[dotsize=0.1,linecolor=RoyalPurple](3.85,4.56)
\psdots[dotsize=0.1,linecolor=RoyalPurple](5.535,4.699)
\psline[linewidth=0.02cm,linecolor=RoyalPurple](2.86875,5.845625)(3.56,6.15)
\psdots[dotsize=0.1,linecolor=MidnightBlue](2.86875,5.845625)
\psdots[dotsize=0.1,linecolor=MidnightBlue](3.26875,4.745625)
\psdots[dotsize=0.1,linecolor=MidnightBlue](4.56875,4.345625)
\psdots[dotsize=0.1,linecolor=MidnightBlue](7.66875,4.145625)
\psdots[dotsize=0.1,linecolor=RoyalPurple](6.96,4.53)
\psdots[dotsize=0.1,linecolor=RoyalPurple](8.26,4.98)
\psdots[dotsize=0.1,linecolor=RoyalPurple](11.465,5.14)
\psdots[dotsize=0.1,linecolor=RoyalPurple](11.73,6.95)
\psdots[dotsize=0.1,linecolor=RoyalPurple](11.6663,7.81)
\psdots[dotsize=0.1,linecolor=RoyalPurple](11.175,8.23)
\psdots[dotsize=0.1,linecolor=RoyalPurple](10.015,8.394375)
\psdots[dotsize=0.1,linecolor=RoyalPurple](8.66,8.494375)
\psdots[dotsize=0.1,linecolor=RoyalPurple](5.03,8.2)
\psdots[dotsize=0.1,linecolor=color82](6.06875,7.445625)
\psdots[dotsize=0.1,linecolor=color82](8.76875,7.245625)
\psdots[dotsize=0.1,linecolor=color82](11.16875,7.245625)
\psdots[dotsize=0.1,linecolor=color82](7.36875,5.645625)
\psdots[dotsize=0.1,linecolor=MidnightBlue](4.06875,7.745625)
\psdots[dotsize=0.1,linecolor=MidnightBlue](5.96875,8.645625)
\psdots[dotsize=0.1,linecolor=MidnightBlue](11.66875,4.445625)
\psdots[dotsize=0.1,linecolor=MidnightBlue](11.86875,7.345625)
\psdots[dotsize=0.1,linecolor=MidnightBlue](11.46875,8.245625)
\psdots[dotsize=0.1,linecolor=MidnightBlue](9.36875,8.945625)
\psline[linewidth=0.04cm](3.06875,3.445625)(2.36875,4.845625)
\psline[linewidth=0.04cm](2.36875,4.845625)(2.46875,6.545625)
\psline[linewidth=0.04cm](2.46875,6.545625)(3.86875,9.045625)
\psline[linewidth=0.04cm](12.46875,6.545625)(11.96875,8.145625)
\psline[linewidth=0.04cm](11.96875,8.145625)(11.16875,9.445625)
\psline[linewidth=0.04cm](8.56875,9.545625)(11.16875,9.445625)
\psline[linewidth=0.04cm](3.86875,9.045625)(8.56875,9.545625)
\psline[linewidth=0.04cm](3.06875,3.445625)(6.56875,3.445625)
\psline[linewidth=0.04cm](6.56875,3.445625)(10.26875,3.445625)
\psline[linewidth=0.04cm](10.26875,3.445625)(11.96875,4.045625)
\psline[linewidth=0.04cm](11.96875,4.045625)(12.46875,6.545625)
\usefont{T1}{ppl}{m}{n}
\rput(2.60375,9.245625){\Huge $\tilde\Delta$}
\end{pspicture}
}

\vspace{-10.5cm}
\caption{Numbering on the vertices in $\tilde\Delta_0^0$.}\label{firstround}
\end{figure}
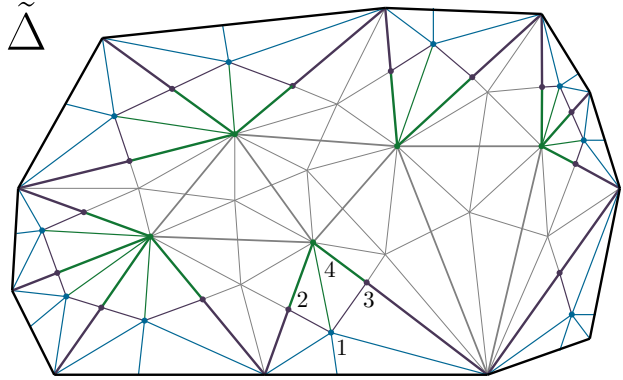

Let us consider the following numbering on the vertices in $\tilde\Delta_0^0$. We take first a triangle $\sigma_i\in\Delta$ with (at least) one of its edges on the boundary, denote by $\gamma_i$ the vertex in $\sigma_i$ that is in the interior of $\Delta$. With the notation as above, we assign the index $1$ to the vertex $\nu_i$, the indices $2$ and $3$ to the vertices $\mu_{ij}$ on the edges of $\sigma_i$, and finally the index $4$ to the vertex $\gamma_i$. After this, we consider the edges of $\sigma_i$ as if they were all part of the boundary and iterate the previous process. We choose another triangle in $\Delta$ having (at least) one edge on the ``new" boundary. We index its vertices following the same order as we did for the ones in $\sigma_i$. We continue this numbering until we have considered all the triangles in $\Delta$ and hence have index all the vertices in $\tilde\Delta^0_0$.  Figure \ref{firstround} shows the triangles that could be considered first, and the different colors refer to the values of $\tilde t_i$ for the interior vertices in each case.

Let us notice that for this numbering, the values for $\tilde t_i$ corresponding to these vertices are either $\geq 3$, as is the case of the vertices $\nu_i$ and $\gamma_i$, or $\tilde t_i=t_i=2$, as it happens for the vertices  $\mu_{ij}$.  Thus, we arrive to the following result.
 \begin{prop}
For any $k\geq 3$:
\[\dim C_k^1(\tilde\Delta)=\binom{k+2}{2}+f^0_1 \binom{k-2}{2}+2f^0_0\binom{k-1}{2}+\binom{2k-1}{2}(f_0-2).\]
For $k=2$: \[\dim C_2^1(\tilde\Delta)=3f_0.\]
\end{prop}
\begin{pf}
It follows directly from the remarks we made above about $\tilde t_i$ for the counting, and Corollary \ref{>r+1}.
\end{pf}

For the numbering we defined on $\tilde\Delta_0^0$,  it is not possible to apply  Theorem \ref{SchTheorem} to find an exact value for the dimension. In fact, when the initial subdivision consists of more than one triangle, by using the numberings allowed in that theorem, we get upper bounds which are strictly bigger that the actual dimension of the space.
\medskip

\noindent{}\textbf{Remark 1.\;}
In this paper we have confined our attention to triangular partitions but it is easy to check that our proofs of the lower and the upper bounds can be extended to rectilinear partitions. The construction also allows to consider regions subdivided by curved boundaries, but then the ideals must be considered individually. The boundary of $|\Delta|$ itself can, of course, be curved.
\medskip

\noindent{}\textbf{Remark 2.\;} The resolution for ideals generated by power of linear forms presented in the paper by \cite{gs}, applied also when the power of the linear forms is different. It makes possible to apply the ideas we present in this paper to mixed splines, i.e., splines where the order of smoothness may differ on the various edges of the subdivision.
\medskip

\noindent{}\textbf{Remark 3.\;} The proof of the exact dimension formula when the degree of the spline space $k$ is $\geq 4r+1$ cannot be directly extended to spaces of polynomials of degree $k\geq 3r+1$. With some restrictions associated to the number of different slopes of the edges containing each vertex it is possible to prove the result for this degree.
\medskip

\noindent{}\textbf{Remark 4.\;} We hope that the two methods we used in the examples to find the exact dimension, namely, by showing that the lower and the upper bounds coincide for certain numbering on the vertices, or that $H_0(\J)=0$ directly by considering the equations of the edges, illustrate the way to easily find the dimension of the spline space for many particular triangulations.
\medskip

\noindent{}\textbf{Remark 5.\;}There is a version of the upper bound
given by \cite{Sch} in \cite{LaiSchu}, where the constraint on the
numbering of the vertices is relaxed. However, for this new version of
the theorem, there are still examples for which the upper bound
formula does not give the correct result, see for instance the example 
p. 243 of \cite{LaiSchu}. 
\medskip 

\noindent{}\textbf{Remark 6.\;} For the three dimensional version of the problem, it is necessary to analyse ideals generated by powers of linear forms in three variables. The  formula of the dimension, analogous to formula (\ref{eq}) in Section \ref{construction}, have two homology modules, and in order to approximate the dimensions the dimension of these two modules must be bound. Some issues and positive results about this problem will be in an upcoming paper.


\end{document}